\documentclass[10pt,letterpaper]{article}
\usepackage{t1enc}
\usepackage[latin1]{inputenc}
\usepackage[english]{babel}
\usepackage{graphics}
\begin{document}
\date{}
\title{   Complete Monotonicity of classical theta functions and applications }
\author{A.  Raouf  Chouikha \footnote
{Universite Paris 13 LAGA UMR 7539 Villetaneuse 93430,  e-mail: chouikha@math.univ-paris13.fr}
}
\maketitle

\begin{abstract}  We produce trigonometric expansions for Jacobi theta functions\\ $\theta_j(u,\tau), j=1,2,3,4$\ where $\tau=i\pi t, t > 0$. This permits us to prove that\ $\log \frac{\theta_j(u, t)}{\theta_j(0, t)}, j=2,3,4$ and $\log \frac{\theta_1(u, t)}{\pi \theta'_1(0, t)}$  as well as $\frac{\frac{\delta\theta_j}{\delta u}}{\theta_j}$ as functions of $t$ are completely monotonic. We also interested in the quotients $S_j(u,v,t) = \frac{\theta_j(u/2,i\pi t)}{\theta_j(u/2,i\pi t)}$. For fixed $u,v$ such that $0\leq u < v < 1$  we prove that the functions $\frac{(\frac{\delta}{\delta t}S_j)}{S_j}$ for $j=1,4$  as well as the functions $-\frac{(\frac{\delta}{\delta t}S_j)}{S_j}$ for $j=2,3$ are completely monotonic for $t \in ]0,\infty[$.\\
{\it Key words and phrases} : theta functions, elliptic functions, complete monotonicity   \footnote {2000 {\it Mathematical Subject Classification} Primary: 33E05, 11F11; Secondary: 33E20, 34A20, 11F20}

\end{abstract}

\section{The theta functions}
Consider the following boundary conditions of the heat equation 
$$\frac{\delta^2 f}{\delta u^2} = \frac{\delta f}{\delta t},\quad f(0,t)= f(1,t), \quad f(u,0)=\pi \delta (u-\frac{1}{2}), \quad 0<u<1$$
where $\delta (u)$ is the Dirac delta function. Then the general solution of the boundary problem is
$$\theta (u,t) = 2 \sum_{n\geq 0} (-1)^n e^{-\frac{{(2n+1)}^2}{2\pi^2 t}} \sin ((2n+1) \pi u)$$

When we write $q = e^{\pi i \tau}$ with $Im \tau > 0$, this solution takes the form 

$$\theta _1(u,\tau ) =\sum_{-\infty \leq n\leq +\infty}(-1)^{\frac {n-1}{2}} q^{(\frac {n+1}{2})^2} e^{i\pi (2n+1)u} = 2 \sum_{n\geq 0} (-1)^n q^{(\frac{n+1}{2})^2} \sin ((2n+1) \pi u)$$
which is the first of the four Jacobi theta functions. It suffices to change the boundary conditions to obtain other theta functions. The corresponding solution of the boundary problem $$\quad \frac{\delta f(0,t)}{\delta u}= \frac{\delta f(0,t)}{\delta u}=0, \quad f(u,0)=\pi \delta (u-\frac{1}{2}), \quad 0<u<1$$ is given by
$$\theta _4(u,\tau ) =\sum_{-\infty \leq n\leq +\infty } (-1)^n q^{n^2} e^{i\pi 2nu} = 1 + 2 \sum_{n\geq 1 } (-1)^n q^{n^2} \cos (2n \pi u)$$
The function $\theta_1(u,\tau))$ is periodic of period 2. We obtain the second theta if we increment $u$ by $1/2$
$$\theta _2(u,\tau ) =\sum_{-\infty \leq n\leq +\infty } q^{(\frac {n+1}{2})^2} e^{i\pi (2n+1)u} = 2 \sum_{n\geq 0} q^{(\frac{n+1}{2})^2} \cos ((2n+1) \pi u)$$
The increment $u$ by $1/2$ yields the third theta function
$$\theta _3(u,\tau ) = 1+2 \sum_{-\infty \leq n\leq +\infty } q^{n^2} \cos (2n \pi u)$$
It is known that these four theta functions can be extended to complex values for $u$ and $q$ such that $\mid q\mid < 1$.\\
All four theta functions are entire functions of $u$. All are periodic, the period of \ $\theta _1$\ and  \ $\theta _2$\ is\ $2$, and that of  \ $\theta _3$\ and  \ $\theta _4$\ is\  $1$. See [1] and [10] for more details\\

It is common knowledge that these functions are fundamental and important and that
they have much extensive applications in different area of mathematical and physical sciences.  Thus, any new property or characterization of theta functions may have consistent implications in fields that use these functions.

\section{Selected facts about completely monotonic functions}
 A real function $f$ is said to be {\it completely monotonic (CM)} on $[0,\infty[$ if $f \ C^\infty$ and all the derivatives verify $(-1)^k f^{(k)}(t) \geq 0, k=0,1;2...$ for all $t > 0$. This definition was introduced in 1921 by F. Hausdorf who called 'total monoton'. The following yields a characterization of such functions, S.N. Bernstein and D. Widder  [3, p.95]

{\it A function $f :[0,\infty[ \rightarrow [0,\infty[$ is CM if and only if there exists a non-decreasing bounded function $h$ such that\  $f(t) = \int_0^\infty e^{-ut} dh(u)$\ and this integral converges for $0<t<\infty$.} \\ $h$ may be considered as a nonegative measure on $[0,\infty[$ such that the integral converges for all $x > 0$.\\ 
Hence, a CM function cannot vanish for any $t > 0$. Notice obviously that the sum and product of CM  functions are CM. Observe that if $f(t)$ is CM then $f^{2m}(t)$ and $-f^{2m+1}(t)$ are also CM. It is known that that if $f(t)$ is CM and $h(t)$ be nonnegative with a CM derivative then $f(h(t))$ also is CM.\\

CM functions have a lot of applications in various fields. These functions have remarkable applications in different branches. In particular, they play a role in potential theory, probability theory, physics, numerical and
asymptotic analysis, and combinatorics, [2].\\

In the recent past, various authors showed that numerous functions, which are defined in terms
of gamma, polygamma, and other special functions, are completely monotonic and used this fact
to derive many interesting new inequalities. Hankel determinant inequality for completely
monotonic functions is proved and it is shown that in connection with an
interpolation problem there exists a close relation between completely monotonic functions and
completely monotonic sequences.
Several related classes of such functions are studied in [2], [12].

There is another class of functions very closed to such functions : the class of logarithmically completely monotonic (LCM). A function $f :[0,\infty[ \rightarrow [0,\infty[$ is logarithmically completely monotonic if it is $C^\infty$ and $(-1)^k (\log f(t))^{(k)} \geq 0, k=1,2,...$. Remark for $k=0$ we may only require that $f(t) > 0.$ Notice that every LCM function is CM. If this inequality  is strict for all $t > 0$  and
$k = 1,2,...$ then f is said to be strictly logarithmically completely monotonic.\\
A function $f$ on $]0,\infty[$ is called a {\it Stieltjes transform} if it can be written in the form
$$f(t) =a + \int_0^{\infty}\frac{d\mu (s)}{s+t}$$
where $a$ is a nonnegative number and $\mu $ a nonnegative measure on $]0,\infty[$ such that \ $\int_0^{\infty}\frac{d\mu (s)}{1+t} < \infty$.\\ It is proved that every Stieltjes transform is a LCM function. It is also proved that if $f$ is LCM then $f^\alpha$ is CM for any $\alpha > 0$.\\
This demonstrates that the investigation of LCM properties of functions are significant and meaningful. 

Various examples of (logarithmically) completely monotonic function have been presented in the litterature; functions associated to Gamma and psi function, a quotients of $K-$Bessel functions and many others, [2].\\

\section{Statements of results}
This paper is organized as follows. At first we descrive expansions of the logarithmic derivatives with respect to $t$ for the four Jacobi theta functions $\theta_j(u,\tau), j=1,2,3,4$ where $\tau=i\pi t, t > 0$. These expressions are derived and inspired from other old ones [4], [5]. We then deduce that these functions  are LCM (logarithmically completely monotonic) for any fixed $u$ such that $0<u<1$ and $t > 0$. We also obtain similar expressions of the logarithmic derivatives with respect to $u$ for the four Jacobi theta functions $\theta_j(u,i\pi t), j=1,2,3,4$ where $t > 0$. Moreover, we prove that $\frac{\frac{d\theta_j}{du}}{\theta_j}, j=1,4$ as functions of $t>0$ are completely monotonic. While $-\frac{\frac{d\theta_j}{du}}{\theta_j}, j=2,3$ are completely monotonic.\\

On the other hand for $u,v \in C$ and $\tau = i\pi t$ with $Re\ t > 0$, define the quotient of theta functions:
$$S_j:=S_j(u,v,t) = \frac{\theta_j(u/2,i\pi t)}{\theta_j(v/2,i\pi t)}.$$
 
Many papers recently interested in monotonicity and convexity of these quotients [6],[7],[10]. This is naturally related to the problem of completely monotonic functions.\\
 A. Solynin and A. Dixit, [7],[10] stated for fixed $u,v$ such that $0\leq u < v < 1$, the functions $S_1(u,v;t)$ and $S_4(u,v;t)$ are positive and strictly increasing for $t \in ]0,\infty[$ while $S_2(u,v;t)$ and $S_3(u,v;t)$ are positive decreasing for $t \in ]0,\infty[$.\\ 
 A. Dixit, A. Roy and A. Zaharescu [7, Th 1.2] proved for $u,v$ such that $0\leq u < v < 1$ the functions $S_2(u,v;t)$ and $S_3(u,v;t)$ are stricly convex for $t \in ]0,\infty[$.\\ 
 
  However, one conjectured [7, Conj 1.1] for $u,v$ such that $0\leq u < v < 1$ the functions $\frac{\delta}{\delta t} S_1(u,v;t), S_2(u,v;t), S_3(u,v;t)$ and $\frac{\delta}{\delta t} S_4(u,v;t)$ are CM for $0 < t < \infty$. In the case where this conjecture is true, as we have seen before by the result of S.N. Bernstein and D. Widder there exist a non-decreasing function $\omega_j$ such that $S_j(u,v;t)=\int_0^\infty e^{-\nu t} d\omega_j(u) d\nu$ for $j = 2, 3$ and $\frac{\delta}{\delta t} S_j(u,v;t)=\int_0^\infty e^{-\nu t} d\omega_j(u) d\nu$ for $j = 3,4$.

Many numerical calculus suggest us that the odd and even derivatives in $t$ of $log(S_j(u,v,t))$ (and not only $S_j(u,v,t)$) have alternating signs. Thus, one naturally may ask if quotients of theta functions are LCM (logarithmically completely monotonic). \\

In this paper we will prove that the quotients $S_j, j=1,2,3,4$ are LCM. More exactly, we prove for $j=2,3$ inequalities $(-1)^k \frac{\delta^k}{\delta t^k}(\frac{\frac{\delta S_j}{\delta t}}{S_j}) < 0$ hold for $k=0,1,2,...$ 
 and for $j=1,4$ inequalities $(-1)^{k} \frac{\delta^k}{\delta t^k}(\frac{\frac{\delta S_j}{\delta t}}{S_j}) >  0$ hold for $k=0,1,2,...$\\ This means  that for fixed $u,v$ such that $0\leq u < v < 1$  the functions $\frac{(\frac{\delta}{\delta t}S_j)}{S_j}$ for $j=1,4$  as well as the functions $-\frac{(\frac{\delta}{\delta t}S_j)}{S_j}$ for $j=2,3$ are completely monotonic for $t \in ]0,\infty[$.
 
Morerover, using properties of theta functions and thank to {\it Maple } or {\it Mathematica} we find again results of [7] and [10]. We prove in addition that for fixed $u,v$ such that $0\leq u < v < 1$ the functions $\frac{\delta}{\delta t} S_4(u,v,t)$ is decreasing and convex and  $\frac{\delta^3}{\delta t^3} S_j(u,v,t), j=2,3$ is non positive for $t \in ]0,\infty[ $.\\

\section{A trigonometric expansions of $\log(\theta_j (u,i\pi t)$}
When $Im \tau >0 $ or equivalently when $t>0$ we proved [3] that the Jacobi theta function $\theta_4$ may be expressed in this trigonometric form which is like a special Fourier expansion
$$ \theta_4(u,\tau) = \theta_4(0,\tau)\ \exp[- \sum_{p\geq 1} \sum_{ k\geq 0} \frac {1}{p} \bigg( \frac {\sin \pi u}{(\sin (k+\frac {1}{2})\pi \tau)}\bigg)^{2p}].$$
By the same way we also obtain similar expansions
$$\theta_3(u,\tau) = \theta_4(0,\tau)\ \exp[- \sum_{p\geq 1} \sum_{ k\geq 0}\frac {1}{p} \bigg( \frac {\cos \pi u}{(\sin (k+\frac {1}{2})\pi \tau)}\bigg)^{2p}],$$
$$\theta_2(u,\tau) = \theta_4(0,\tau)\ \exp[ i\pi (u+\frac {1}{4}\tau)- \sum_{p\geq 1} \sum_{ k\geq 0}\frac {1}{p} \bigg( \frac {\cos \pi (u+\frac {1}{2}\tau)}{(\sin (k+\frac {1}{2})\pi \tau)}\bigg)^{2p}],$$
$$\theta_1(u,\tau) = \theta_4(0,\tau)\ \exp[ i\pi (u-\frac {1}{2}+\frac {1}{4}\tau)- \sum_{p\geq 1} \sum_{k\geq 0}\frac {1}{p} \bigg( \frac {\sin \pi (u+\frac {1}{2}\tau)}{(\sin (k+\frac {1}{2})\pi \tau)}\bigg)^{2p}].$$ 
The above expressions of \ $\theta_3 $ \ and \ $ \theta_4$ \ are valid in the "strip" \ $\mid Im u \mid < \frac{1}{2} Im \tau ,$\\ those relating to  \ $\theta_1 $ \ and \ $ \theta_2$ \  are valid in the "strip" \ $\mid Im u \mid < Im  \tau .$ \\

By the change $\tau = i \pi t$ these expressions may be rewritten as functions of $t >0$. We then state the following\\ 

{\bf Proposition 1}\ {\it When  $t>0$  the Jacobi theta functions $\theta_j, j=1,2,3,4$ may be expressed in this trigonometric form 
$$ \theta_4(u,t) = \theta_4(0,t)\ \exp[- \sum_{p\geq 1} \sum_{ k\geq 0} \frac {(-1)^p}{p} \bigg( \frac {\sin \pi u}{(\sinh (k+\frac {1}{2})\pi^2 t)}\bigg)^{2p}],$$
$$\theta_3(u,t) = \theta_4(0,t)\ \exp[- \sum_{p\geq 1} \sum_{ k\geq 0}\frac {(-1)^p}{p} \bigg( \frac {\cos \pi u}{(\sinh (k+\frac {1}{2})\pi^2 t)}\bigg)^{2p}],$$
$$\theta_2(u,t) = \theta_4(0,t)\ \exp[ i\pi (u-\frac {\pi t}{4})- \sum_{p\geq 1} \sum_{ k\geq 0}\frac {(-1)^p}{p} \bigg( \frac {\cosh \pi (i u-\frac {\pi}{2}t)}{(\sinh (k+\frac {1}{2})\pi^2)}\bigg)^{2p}],$$
$$\theta_1(u,t) = \theta_4(0,t)\ \exp[ i\pi (u-\frac {1}{2}-\frac {i\pi t}{4})- \sum_{p\geq 1} \sum_{k\geq 0}\frac {1}{p} \bigg( \frac {\sinh \pi (i u-\frac {\pi t}{2})}{(\sinh (k+\frac {1}{2})\pi^2 t)}\bigg)^{2p}].$$ 
The above expressions of \ $\theta_3 $ \ and \ $ \theta_4$ \ are valid in the "strip" \ $\mid Im u \mid < \frac{1}{2} Im \tau ,$\\ those relating to  \ $\theta_1 $ \ and \ $ \theta_2$ \  are valid in the "strip" \ $\mid Im u \mid < Im  \tau .$ }\\

Using the same technics as [3],[4] many others similar expansions may be found. In particular, we may prove the following\\

{\bf Proposition 2}\ {\it When  $t>0$  the Jacobi theta functions $\theta_j, j=1,2,3,4$ may be expressed in this trigonometric form 
$$ \theta_2(u,t) = \theta_2(0,t)\ \exp[\log{\cos \pi u}- \sum_{p\geq 1} \sum_{ k\geq 1} \frac {(-1)^p}{p} \bigg( \frac {\sin \pi u}{(\sinh (k\pi^2 t)}\bigg)^{2p}],$$
$$\theta_1(u,t) = \pi \frac{\delta \theta_1}{\delta u}(0,t)\ \exp[\log{\sin \pi u}- \sum_{p\geq 1} \sum_{ k\geq 1} \frac {(-1)^p}{p} \bigg( \frac {\cos \pi u}{(\sinh (k\pi^2 t)}\bigg)^{2p}],$$
$$\theta_3(u,t) = \theta_2(0,t) \exp [(i\pi u+i\pi^2 \frac{t}{4}+ \log{\cos (\pi u+i\pi^2 \frac{t}{2})}- \sum_{p\geq 1} \sum_{ k\geq 1} \frac {(-1)^p}{p} \bigg( \frac {\sin (\pi u+i\pi^2 \frac{t}{2})}{(\sinh (k\pi^2 t)}\bigg)^{2p}],$$
$$\theta_4(u,t) = \theta_2(0,t) \exp [(i\pi u+i\pi^2 \frac{t}{4}+ \log{\sin (\pi u+i\pi^2 \frac{t}{2})}- \sum_{p\geq 1} \sum_{ k\geq 1} \frac {(-1)^p}{p} \bigg( \frac {\cos (\pi u+i\pi^2 \frac{t}{2})}{(\sinh (k\pi^2 t)}\bigg)^{2p}],$$
 above expressions of \ $\theta_3 $ \ and \ $ \theta_4$ \ are valid in the "strip" \ $\mid Im u \mid < \frac{1}{2} Im \tau ,$\\ those relating to  \ $\theta_1 $ \ and \ $ \theta_2$ \  are valid in the "strip" \ $\mid Im u \mid < Im  \tau .$ }\\ 

Using the above expansions turn out now to express logarithmic derivatives with respect to $t > 0$ of the theta functions. In the sequel we will denote $\theta_j(u,\tau)=\theta_j(u,i \pi t)=\theta_j(u,t)$. At first, we prove the following\\

{\bf Theorem 1}\ {\it The logarithmic derivatives with respect to $t > 0$ of the Jacobi theta functions $\frac{\theta'_j(u, t)}{\theta_j(u, t)}$  where $\theta'_j(u, t)=\frac{\delta \theta_{j}(u, t)}{\delta t}$ may be expressed under the following form
$$-\frac{\theta'_3}{\theta_3}(u,t)+\frac{\theta'_4}{\theta_4}(0,t)= \sum_{k\geq 0}\frac{\pi^2 ( 2k+1 ) [\coth ( k
+\frac{1}{2} ) \pi^2 \,t ][\cos \pi {u}]^{2}}{ [(\sinh \left(  \left( k
+\frac{1}{2} \right) \pi^2 \,t \right)]^2+[\cos \pi {u}]^{2}}=2\sum_{k\geq 0}\frac{\pi^2 ( 2k+1 ) [\coth ( k
+\frac{1}{2} ) \pi^2 \,t ][\cos \pi {u}]^{2}}{ \left(\cosh \left(  \left( 2k+1
 \right) \pi^2 \,t \right)+[\cos 2\pi {u}]\right)}, $$
$$-\frac{\theta'_4}{\theta_4}(u,t) +\frac{\theta'_4}{\theta_4}(0,t)= \sum_{k\geq 0}\frac{\pi^2 ( 2k+1 ) [\coth ( k
+\frac{1}{2} ) \pi^2 \,t ][\sin \pi u]^{2}}{ [\sinh \left(  \left( k
+\frac{1}{2} \right) \pi^2 \,t \right)]^2+[\sin \pi {u}]^{2}}= 2\sum_{k\geq 0}\frac{\pi^2 ( 2k+1 ) [\coth ( k
+\frac{1}{2} ) \pi^2 \,t ][\sin \pi u]^{2}}{ \left(\cosh \left(  \left( 2k+1
 \right) \pi^2 \,t \right)-[\cos 2\pi {u}]\right)}$$
where $\mid Im u \mid < \frac{1}{2} \pi t ,$
$$-\frac{\theta'_1}{\theta_1}(u,t)+{\pi \frac{{\theta'_1}}{\frac{\delta \theta_1}{\delta u}}}(0,t)=\sum_{k\geq 1}\frac{\pi^2 ( 2k ) [\coth ( k
 \pi^2 \,t )][\sin \pi u]^{2}}{ [(\sinh   \left( k
  \pi^2 \,t \right) ]^2+[\sin \pi {u}]^{2}}= 2\sum_{k\geq 1}\frac{\pi^2 ( 2k ) [\coth ( k
 \pi^2 \,t )][\sin \pi u]^{2}}{ \left(\cosh \left(  \left( 2k
 \right) \pi^2 \,t \right)-[\cos 2\pi {u}]\right)}$$
$$-\frac{\theta'_2}{\theta_2}(u,t)+\frac{\theta'_2}{\theta_2}(0,t)=\sum_{k\geq 1}\frac{\pi^2 ( 2k ) [\coth ( k
  \pi^2 \,t )][\cos \pi {u}]^{2}}{ [(\sinh   \left( k
  \pi^2 \,t \right) ]^2+[\cos \pi {u}]^{2}}=2\sum_{k\geq 1}\frac{\pi^2 ( 2k ) [\coth ( k
  \pi^2 \,t )][\cos \pi {u}]^{2}}{ \left(\cosh \left(  \left( 2k
 \right) \pi^2 \,t \right)+[\cos 2\pi {u}]\right)}$$
where} \ $\mid Im u \mid <  \pi t $.\\

{\bf Proof}\\
It suffices to prove expansions for $\theta_4$ and $\theta_2$ for example. The others will easily be deduced from the one choosen. Consider by Proposition 1   
$$\frac{\theta_4(u, t)}{\theta_4(0, t)} = exp-[\sum_{k\geq 0}\sum_{p\geq 1}\frac {(-1)^p}{p} \frac {[\sin \pi u]^{2p}}{[\sinh (k+\frac {1}{2})\pi^2 t]^{2p}}]$$ where we replaced $\tau$ by $i\pi t$. We then have

$$\frac{\frac{\delta \theta_4(u, t)}{\delta t}}{\theta_4(u, t)}-\frac{\theta'_4}{\theta_4}(0,t)= \sum_{k\geq 0}\sum_{p\geq 1} {\frac { \left( -1 \right) ^{p} \cosh \left(  \left( k+\frac{1}{2} \right) \pi^2 \,t \right)  \left( 2\,k+1 \right) \pi^2 }{\left( \sinh \left(  \left( k
+\frac{1}{2} \right) \pi^2 \,t \right)  \right) ^{2\,p+1}}}(-[\sin \pi u]^{2p})$$

$$=\sum_{k\geq 0}\sum_{p\geq 1}  \pi^2 \left( -1 \right) ^{p}\left( 2k+1 \right) \left( \sinh \left(  \left( k
+\frac{1}{2} \right) \pi^2 \,t \right)  \right) ^{-2\,p-1} [\cosh \left( k
+\frac{1}{2} \right) \pi^2 \,t ] (- [\sin \pi u]^{2p})  $$
$$= \pi^2 \left( 2k+1 \right) [\cosh \left( k
+\frac{1}{2} \right) \pi^2 \,t ]  \sum_{p\geq 1}  \left( -1 \right) ^{p} (-[\sin \pi u]^{2p}) \left( \sinh \left(  \left( k
+\frac{1}{2} \right) \pi^2 \,t \right)  \right) ^{-2\,p-1}$$
Compute now the sum 
$$A_{k,4} =\sum_{p\geq 1}  \left( -1 \right) ^{p} (-[\sin \pi u]^{2p}) \left( \sinh \left(  \left( k
+\frac{1}{2} \right) \pi^2 \,t \right)  \right) ^{-2\,p-1}=$$ 
$$-\left( \sinh \left(  \left( k
+\frac{1}{2} \right) \pi^2 \,t \right)^{-1}\right)\sum_{p\geq 1}  \left( -1 \right) ^{p}{\left( \sinh \left(  \left( k
+\frac{1}{2} \right) \pi^2 \,t \right)  \right) ^{2\,p}} \frac{[\sin \pi u]^{2p})}{\left( \sinh \left(  \left(k
+\frac{1}{2} \right) \pi^2 \,t \right)  \right) ^{2\,p}}]= $$

$$\frac{[\sin \pi u]^{2}}{\left( \sinh \left(  \left( k
+\frac{1}{2} \right) \pi^2 \,t \right)  \right) ^{3}} \sum_{p\geq 0} \frac{(-1)^{p}[\sin \pi u]^{2p}}{\left( \sinh \left(  \left( k
+\frac{1}{2} \right) \pi^2 \,t \right)  \right) ^{2\,p}}$$
Since the variables $u$ belong to the "strip" \ $\mid Im u \mid < \frac{1}{2} \pi t ,$ it implies that $\frac{[\sin \pi u]}{\left( \sinh \left(  \left( k
+\frac{1}{2} \right) \pi^2 \,t \right)  \right) } < 1.$ Then, the sum equals
$$A_{k,4} = - \frac{[\sin \pi u]^{2}}{\left( \sinh \left(  \left( k
+\frac{1}{2} \right) \pi^2 \,t \right)  \right) ^{3}} \frac{1}{1+\frac{[\sin \pi u]^{2}}{\left( \sinh \left(  \left( k
+\frac{1}{2} \right) \pi^2 \,t \right)  \right) ^{2}}}$$
After simplifying we get
$$A_{k,4}= -\frac{\left( \sinh \left(  \left( k
+\frac{1}{2} \right) \pi^2 \,t \right)  \right)^{-1}[\sin \pi u]^{2}}{ [\sinh \left(  \left( k
+\frac{1}{2} \right) \pi^2 \,t \right)]^2+[\sin \pi u]^{2}} $$
So we obtain the expression for the derivative $$\frac{\theta'_4}{\theta_4}(u,t)-\frac{\theta'_4}{\theta_4}(0,t)= -\sum_{k\geq 0} \pi^2 ( 2k+1 ) [\coth ( k
+\frac{1}{2} ) \pi^2 \,t ]A_{k,4}$$
 
$$\frac{\theta'_4}{\theta_4}(u,t)-\frac{\theta'_4}{\theta_4}(0,t) = -\sum_{k\geq 0}\frac{\pi^2 ( 2k+1 ) [\coth ( k
+\frac{1}{2} ) \pi^2 \,t ][\sin \pi u]^{2}}{ [\sinh \left(  \left( k
+\frac{1}{2} \right) \pi^2 \,t \right)]^2+[\sin \pi {u}]^{2}}$$ We then easily deduce the expansion for $$\frac{\theta'_3}{\theta_3}(u,t)-\frac{\theta'_4}{\theta_4}(0,t)= -\sum_{k\geq 0}\frac{\pi^2 ( 2k+1 ) [\coth ( k
+\frac{1}{2} ) \pi^2 \,t ][\cos \pi {u}]^{2}}{ [(\sinh \left(  \left( k
+\frac{1}{2} \right) \pi^2 \,t \right)]^2+[\cos \pi {u}]^{2}}$$ since $\frac{\theta'_4}{\theta_4}(u+\frac{1}{2},t) =\frac{\theta'_3}{\theta_3}(u,t).$\\
Consider now by Proposition 2  
$$\theta_2(u,t) = \theta_2(0,t)\ \exp[\log{\cos \pi u}- \sum_{p\geq 1} \sum_{ k\geq 1} \frac {(-1)^p}{p} \bigg( \frac {\sin \pi u}{(\sinh (k\pi^2 t)}\bigg)^{2p}].$$ We then have 
$$\frac{\frac{\delta \theta_2(u, t)}{\delta t}}{\theta_2(u, t)}-\frac{\theta'_2}{\theta_2}(0,t)=  \sum_{k\geq 0}\sum_{p\geq 1} {\frac { \left( -1 \right) ^{p} \cosh \left(  \left( k \right) \pi^2 \,t \right)  \left( 2\,k \right) \pi^2 }{\left( \sinh \left(  \left( k
 \right) \pi^2 \,t \right)  \right) ^{2\,p+1}}}(-[\sin \pi u]^{2p})$$
 
 $$=\sum_{k\geq 0}\sum_{p\geq 1}  \pi^2 \left( -1 \right) ^{p}\left( 2k \right) \left( \sinh \left(  \left( k
 \right) \pi^2 \,t \right)  \right) ^{-2\,p-1} [\cosh \left( k
 \right) \pi^2 \,t ] (- [\sin \pi u]^{2p})  $$
$$= \pi^2 \left( 2k \right) [\cosh \left( k
 \right) \pi^2 \,t ]  \sum_{p\geq 1}  \left( -1 \right) ^{p} (-[\sin \pi u]^{2p}) \left( \sinh \left(  \left( k
 \right) \pi^2 \,t \right)  \right) ^{-2\,p-1}$$
Compute now the sum 
$$B_{k,4} =\sum_{p\geq 1}  \left( -1 \right) ^{p} (-[\sin \pi u]^{2p}) \left( \sinh \left(  \left( k
 \right) \pi^2 \,t \right)  \right) ^{-2\,p-1}=$$ 
$$-\left( \sinh \left(  \left( k
 \right) \pi^2 \,t \right)^{-1}\right)\sum_{p\geq 1}  \left( -1 \right) ^{p}{\left( \sinh \left(  \left( k
 \right) \pi^2 \,t \right)  \right) ^{2\,p}} \frac{[\sin \pi u]^{2p})}{\left( \sinh \left(  \left(k
 \right) \pi^2 \,t \right)  \right) ^{2\,p}}]= $$
$$\frac{[\sin \pi u]^{2}}{\left( \sinh \left(  \left( k
 \right) \pi^2 \,t \right)  \right) ^{3}} \sum_{p\geq 0} \frac{(-1)^{p}[\sin \pi u]^{2p}}{\left( \sinh \left(  \left( k
+\frac{1}{2} \right) \pi^2 \,t \right)  \right) ^{2\,p}}$$
Since the variables $u$ belong to the "strip" \ $\mid Im u \mid <  \pi t ,$ it implies that $\frac{[\sin \pi u]}{\left( \sinh \left(  \left( k
 \right) \pi^2 \,t \right)  \right) } < 1.$ Then, the sum equals
$$B_{k,4} = - \frac{[\sin \pi u]^{2}}{\left( \sinh \left(  \left( k
 \right) \pi^2 \,t \right)  \right) ^{3}} \frac{1}{1+\frac{[\sin \pi u]^{2}}{\left( \sinh \left(  \left( k
 \right) \pi^2 \,t \right)  \right) ^{2}}}$$
After simplifying we get
$$B_{k,4}= -\frac{\left( \sinh \left(  \left( k
 \right) \pi^2 \,t \right)  \right)^{-1}[\sin \pi u]^{2}}{ [\sinh \left(  \left( k
 \right) \pi^2 \,t \right)]^2+[\sin \pi u]^{2}} $$
So we obtain the expression for the derivative $$\frac{\theta'_2}{\theta_2}(u,t)-\frac{\theta'_4}{\theta_4}(0,t)= -\sum_{k\geq 0} \pi^2 ( 2k ) [\coth ( k
 ) \pi^2 \,t ]B_{k,4}$$
 
$$\frac{\theta'_2}{\theta_2}(u,t)-\frac{\theta'_2}{\theta_2}(0,t) = -\sum_{k\geq 0}\frac{\pi^2 ( 2k ) [\coth ( k
 ) \pi^2 \,t ][\sin \pi u]^{2}}{ [\sinh \left(  \left( k
 \right) \pi^2 \,t \right)]^2+[\sin \pi {u}]^{2}}$$ We then easily deduce the expansion for $$\frac{\theta'_1}{\theta_1}(u,t)-\frac{\theta'_2}{\theta_2}(0,t)= -\sum_{k\geq 0}\frac{\pi^2 ( 2k ) [\coth ( k
 ) \pi^2 \,t ][\cos \pi {u}]^{2}}{ [(\sinh \left(  \left( k
 \right) \pi^2 \,t \right)]^2+[\cos \pi {u}]^{2}}$$ since $\frac{\theta'_2}{\theta_2}(u+\frac{1}{2},t) =\frac{\theta'_1}{\theta_1}(u,t).$\\
 
\bigskip

 {\bf Theorem 2}\ {\it The logarithmic derivatives with respect to  $u$ of the Jacobi theta functions $\frac{\frac{\delta \theta_{j}(u, t)}{\delta u}}{\theta_j(u, t)}$ may be also expressed 
$$ \frac{\frac{\delta \theta_{4}}{\delta u}}{\theta_4}(u, t)=2\sum_{k\geq 0}\,{\frac {\pi \,\sin \left(2 \pi \,u \right) }{\cosh \left( 2\,\pi^2 \,t
k+\pi^2 \,t \right) -\cos \left( \pi \,u \right) }}$$
$$ \frac{\frac{\delta \theta_{3}}{\delta u}}{\theta_3}(u, t)=-2\sum_{k\geq 0}\,{\frac {\pi \,\sin \left(2 \pi \,u \right) }{\cosh \left( 2\,\pi^2 \,t
k+\pi^2 \,t \right) +\cos \left( \pi \,u \right) }}.$$
$$ \frac{\frac{\delta \theta_{1}}{\delta u}}{\theta_1}(u, t)=\cot \pi u+2\pi \sum_{k\geq 1}\frac{\sin (2\pi u)}{\cosh (2k\pi^2 t)-\cos (2\pi u)}$$
$$ \frac{\frac{\delta \theta_{2}}{\delta u}}{\theta_2}(u, t)=-\tan \pi u-2\pi \sum_{k\geq 1}\frac{\sin (2\pi u)}{\cosh (2k\pi^2 t)+\cos (2\pi u)}$$
where $t>0$ and $0<u<1$}.\\

{\bf Proof}\\
It suffices to prove expansions for $\theta_4$ and $\theta_2$ for example. The others will be easily deduced. Consider by Proposition 1  
$$\frac{\theta_4(u, t)}{\theta_4(0, t)} = exp-[\sum_{k\geq 0}\sum_{p\geq 1}\frac {(-1)^p}{p} \frac {[\sin \pi u]^{2p}}{[\sinh (k+\frac {1}{2})\pi^2 t]^{2p}}]$$. We then have

$$\frac{\frac{\delta \theta_4(u, t)}{\delta u}}{\theta_4(u, t)}=-\sum_{k\geq 0}\sum_{p>0}{\frac { \left( -1 \right) ^{p} \left( \sin \left( \frac {\pi}{2} \,u
 \right)  \right) ^{2\,p}\cos \left( \frac {\pi}{2} \,u \right) \pi^2 }{\sin
 \left( \frac {\pi}{2} \,u \right)  \left( \sinh \left( \pi^2 \, \left( k+\frac {1}{2}
 \right) t \right)  \right) ^{2\,p}}}. $$
 Since the variables $u$ belong to the "strip" \ $\mid Im u \mid < \frac{1}{2} \pi t ,$ it implies that $\frac{[\sin \pi u]}{\left( \sinh \left(  \left( k
+\frac{1}{2} \right) \pi^2 \,t \right)  \right) } < 1.$ Then, the sum equals
 $$ \sum_{k\geq 0}\left( -1 \right) ^{p} \left( {\frac {\sin \left( \frac {\pi}{2} \,u
 \right) }{\sinh \left( \pi^2 \, \left( k+\frac {1}{2} \right) t \right) }}
 \right) ^{2\,p}\cos \left( \frac {\pi}{2} \,u \right) \pi^2  \left( \sin
 \left( \frac {\pi}{2} \,u \right)  \right) ^{-1}=$$
 $$\sum_{k\geq 0}\pi^2 \,\sin \left( \frac {\pi}{2} \,u \right) \cos \left( \frac {\pi}{2} \,u
 \right)  \left( \sinh \left( \pi^2 \, \left( k+\frac {1}{2} \right) t \right) 
 \right) ^{-2} \left( 1+{\frac { \left( \sin \left( \frac {\pi}{2} \,u
 \right)  \right) ^{2}}{ \left( \sinh \left( \pi^2 \, \left( k+\frac {1}{2}
 \right) t \right)  \right) ^{2}}} \right) ^{-1}.$$
 After simplifying we get
$$\frac{\frac{\delta \theta_4}{\delta u}}{\theta_4}(u,t)=\sum_{k\geq 0}{\frac {2\pi^2 \,\sin \left( \frac {\pi}{2} \,u \right) \cos \left( \frac {\pi}{2} \,u \right) }{ \left( \cosh \left( \pi^2 \,tk+\frac {1}{2}\,\pi^2 \,t \right) 
 \right) ^{2}- \left( \cos \left( \frac {\pi}{2} \,u \right)  \right) ^{2}}}=2\sum_{k\geq 0}\,{\frac {\pi^2 \,\sin \left( \pi \,u \right) }{\cosh \left( 2\,\pi^2 \,t
k+\pi \,t \right) -\cos \left( \pi \,u \right) }}.$$
We then deduce the expansion for ${\theta_3}$
$$ \frac{\frac{\delta \theta_{3}}{\delta u}}{\theta_3}(u, t)=-2\sum_{k\geq 0}\,{\frac {\pi \,\sin \left(2 \pi \,u \right) }{\cosh \left( 2\,\pi^2 \,t
k+\pi^2 \,t \right) +\cos \left( \pi \,u \right) }}.$$

Consider now by Proposition 2
$$\frac{\theta_2(u, t)}{\theta_2(0, t)} = exp-[\sum_{k\geq 0}\sum_{p\geq 1}\frac {(-1)^p}{p} \frac {[\sin \pi u]^{2p}}{[\sinh (k)\pi^2 t]^{2p}}]$$. We then have

$$\frac{\frac{\delta \theta_2(u, t)}{\delta u}}{\theta_2(u, t)}=-\sum_{k\geq 0}\sum_{p>0}{\frac { \left( -1 \right) ^{p} \left( \sin \left( \frac {\pi}{2} \,u
 \right)  \right) ^{2\,p}\cos \left( \frac {\pi}{2} \,u \right) \pi^2 }{\sin
 \left( \frac {\pi}{2} \,u \right)  \left( \sinh \left( \pi^2 \, \left( k
 \right) t \right)  \right) ^{2\,p}}}. $$
 Since the variables $u$ belong to the "strip" \ $\mid Im u \mid <  \pi t ,$ it implies that $\frac{[\sin \pi u]}{\left( \sinh \left(  \left( k
 \right) \pi^2 \,t \right)  \right) } < 1.$ Then, the sum equals
 $$ \sum_{k\geq 0}\left( -1 \right) ^{p} \left( {\frac {\sin \left(  \,u
 \right) }{\sinh \left( \pi^2 \, \left( k \right) t \right) }}
 \right) ^{2\,p}\cos \left( \frac {\pi}{2} \,u \right) \pi^2  \left( \sin
 \left( \frac {\pi}{2} \,u \right)  \right) ^{-1}=$$
 $$\sum_{k\geq 0}\pi^2 \,\sin \left( \frac {\pi}{2} \,u \right) \cos \left( \frac {\pi}{2} \,u
 \right)  \left( \sinh \left( \pi^2 \, \left( k \right) t \right) 
 \right) ^{-2} \left( 1+{\frac { \left( \sin \left( \frac {\pi}{2} \,u
 \right)  \right) ^{2}}{ \left( \sinh \left( \pi^2 \, \left( k
 \right) t \right)  \right) ^{2}}} \right) ^{-1}.$$
 After simplifying we get
$$ \frac{\frac{\delta \theta_{2}}{\delta u}}{\theta_2}(u, t)=\sum_{k\geq 0}{\frac {2\pi^2 \,\sin \left( \frac {\pi}{2} \,u \right) \cos \left( \frac {\pi}{2} \,u \right) }{ \left( \cosh \left( \pi^2 \,tk\,\pi^2 \,t \right) 
 \right) ^{2}- \left( \cos \left( \frac {\pi}{2} \,u \right)  \right) ^{2}}}=2\sum_{k\geq 0}\,{\frac {\pi^2 \,\sin \left( \pi \,u \right) }{\cosh \left( 2\,\pi^2 \,t
k+\pi \,t \right) -\cos \left( \pi \,u \right) }}.$$
We then deduce the expansion for $\frac{\theta'_1}{\theta_1}.$
$$ \frac{\frac{\delta \theta_{1}}{\delta u}}{\theta_1}(u, t)=\cot \pi u+2\pi \sum_{k\geq 1}\frac{\sin (2\pi u)}{\cosh (2k\pi^2 t)-\cos (2\pi u)}$$

\bigskip
{\bf Remark}\quad The logarithmic derivative of theta functions is related to the Zeta functions of Jacobi. Indeed, it is defined by $$Zn(z,k) = \frac{1}{2k}\frac{\delta}{\delta z}\log \theta_4(u,\tau)$$ where $u = \frac{z}{2K}$ and $K=2\int_0^{\frac{\pi}{2}}\frac{dx}{\sqrt{1-k^2 \sin^2x}}$ is the complete elliptic integral of the first kind and the modulus is such that $0 < k < 1$.\\
It is wellknown that the Zeta function has a Fourier expansion 
$$Zn(z,k) = \frac{2\pi}{K}\sum_{n\geq 1 } \frac{q^n}{1-q^{2n}}\sin(\frac{n\pi z}{K})$$ which may be rewritten by Theorem 2
$$Zn(z,k) = \frac{2\pi}{K}\sum_{k\geq 0}\,{\frac {\pi^2 \,\sin \left( \pi \,\frac{z}{2K} \right) }{\cosh \left( 2\,\pi^2 \,t
k+\pi^2 \,t \right) -\cos \left( \pi \,\frac{z}{2K} \right) }}.$$\\

Our first main result is the following
 
\bigskip
 
{\bf Theorem 3}\quad {\it For any fixed $u$ such that $0<u<1$ the Jacobi theta functions $\frac{\theta_j(u, t)}{\theta_j(0, t)}, j=2,3,4$ and $\frac{\theta_1(u, t)}{\pi \theta'_1(0, t)}$ (where $\theta'_1(0, t) = \frac{\delta \theta_1}{\delta u}(O,t)$) are logarithmically completely monotonic with respect to $t > 0$}\\

\bigskip

This theorem may be deduced from the following 
\bigskip

{\bf Lemma 1}\quad {\it The following functions are LCM (logarithmically completely monotonic) for } $t>0$:\\
 $$ \coth (t),\quad  \frac{1}{\sinh (t)}, \quad \frac{a}{\cosh (t)-b}, a>0,\ 0\leq b\leq 1$$
 
 \bigskip
 {\bf Proof of Lemma 1}\quad We may verify that these functions are CM. To insure that they are LCM we follow  Chun-Fu Wei1, Bai-Ni Guo [9] who considered the ith derivative of the functions $\frac{1}{e^{\mp t} - 1}$: $$F_i(t) = (-1)^i\frac{d^i}{dt^i}(\frac{1}{e^t - 1}) \ and \ G_i(t) = (-1)^i\frac{d^i}{dt^i}(\frac{1}{1-e^{-t}})$$ where \ $i=0,1,2,...$. They proved that $F_i, G_i$  are CM for all $t > 0$. But  $$F_0(t) = \frac{1}{e^t - 1}\ and \ G_0(t) = \frac{1}{1 - e^{-t}}$$ are LCM for all $t > 0$. Consequently, $$1 +2 F_0(t) = \frac{e^t +1}{e^t - 1}= \coth \frac{t}{2}$$ is also LCM for all $t > 0$ as well as $$\frac{1}{\sinh \frac{t}{2}}= 2 e^t F_0(t) = 2 e^{-t} G_0(t)\ and \ \frac{a}{\cosh (t)-b}$$ are also LCM.\\ 
  
This Lemma as well as Theorem 1 implies Theorem 3.

\bigskip
By the same way and using again this Lemma we prove the following

\bigskip
 
 {\bf Theorem 4}\quad {\it For any fixed $u$ such that $0<u<1$ the derivatives of Jacobi theta functions $\frac{\frac{\delta \theta_{j}(u, t)}{\delta u}}{\theta_j(u, t)}$  are completely monotonic with respect to $t > 0$ for $j=1,4$. While $-\frac{\frac{\delta \theta_{j}(u, t)}{\delta u}}{\theta_j(u, t)}$ are completely monotonic with respect to $t > 0$ for $j=2,3$}

 \bigskip
 
 {\bf Proof} \ By Theorem 2 we have $ \frac{\frac{\delta \theta_{4}}{\delta u}}{\theta_4}(u, t)=2\sum_{k\geq 0}\,{\frac {\pi \,\sin \left(2 \pi \,u \right) }{\cosh \left( 2\,\pi^2 \,t
k+\pi^2 \,t \right) -\cos \left( \pi \,u \right) }}$. It implies by the preceding Lemma 1 that $\frac{\frac{\delta \theta_{4}}{\delta u}}{\theta_4}(u, t)$ is CM with respect to $t > 0$. It is the same for $\frac{\frac{\delta \theta_{1}}{\delta u}}{\theta_1}(u, t), -\frac{\frac{\delta \theta_{2}}{\delta u}}{\theta_2}(u, t)$ and $-\frac{\frac{\delta \theta_{3}}{\delta u}}{\theta_3}(u, t).$\\

\section {Some applications for the quotients of theta functions}
For $u,v \in C$ and $\tau = i\pi t$ with $Re t > 0$, we define the quotient of theta functions as follows
$$S_j:=S_j(u,v,t) = \frac{\theta_j(u/2,i\pi t)}{\theta_j(v/2,i\pi t)}.$$
Many papers recently interested in monotonicity and convexity of these quotients [6],[7],[10]. This is related to the problem of completely monotonic functions.\\
 A. Solynin and A. Dixit, [7],[10] proved the monotonicity of $S_j(u,v,t) = \frac{\theta_j(u/2,i\pi t)}{\theta_j(v/2,i\pi t)}$. More precisely they stated for fixed $u,v$ such that $0\leq u < v < 1$, the functions $S_1(u,v;t)$ and $S_4(u,v;t)$ are positive and strictly increasing for $t \in ]0,\infty[$ while $S_2(u,v;t)$ and $S_3(u,v;t)$ are positive decreasing for $t \in ]0,\infty[$.\\ 
 A. Dixit, A. Roy and A. Zaharescu [7, Th 1.2] proved for $u,v$ such that $0\leq u < v < 1$ the functions $S_2(u,v;t)$ and $S_3(u,v;t)$ are stricly convex for $t \in ]0,\infty[$.\\ 
  However, one conjectured [7, Conj 1.1] for $u,v$ such that $0\leq u < v < 1$ the functions $\frac{\delta}{\delta t} S_1(u,v;t), S_2(u,v;t), S_3(u,v;t)$ and $\frac{\delta}{\delta t} S_4(u,v;t)$ are CM for $0 < t < \infty$.\\
    
Many numerical calculus suggest us that the odd and even derivatives in $t$ of $log(S_j(u,v,t))$ have alternating signs. Thus, one naturally ask if quotients of theta functions are LCM (logarithmically completely monotonic).\\
  Using the above expansions of theta functions we shall solve this problem. Moreover, as we have seen before by the result of S.N. Bernstein and D. Widder there exists a non-decreasing function $\omega_j$ such that $S_j(u,v;t)=\int_0^\infty e^{-\nu t} d\omega_j(u) d\nu$ for $j = 2, 3$ and $\frac{\delta}{\delta t} S_j(u,v;t)=\int_0^\infty e^{-\nu t} d\omega_j(u) d\nu$ for $j = 3,4$.\\ 
 
\bigskip
 
 {\bf Theorem 5}\quad {\it For fixed $u,v$ such that $0\leq u < v < 1$, the quotients of theta functions $S_j:=S_j(u,v,t) = \frac{\theta_j(u/2, t)}{\theta_j(u/2, t)}, j=1,2,3,4$ are strictly LCM (logarithmically completely monotonic) for $t > 0.$ That means :\\ 
 - for $j=2,3$ inequalities $(-1)^k \frac{\delta^k}{\delta t^k}(\frac{\frac{\delta S_j}{\delta t}}{S_j}) < 0$ hold for $k=0,1,2,...$\\ 
 - for $j=1,4$ inequalities $(-1)^{k} \frac{\delta^k}{\delta t^k}(\frac{\frac{\delta S_j}{\delta t}}{S_j}) >  0$ hold for $k=0,1,2,...$ }\\
 
 \bigskip
 To prove Theorem 5 we need this lemma\\
 
 {\bf Lemma 2}\quad {\it Let $\theta_j(u,i\pi t), j=1,2,3,4$ the four theta fonctions. then the following identities
 $$\frac{\delta }{\delta u}(\frac{\delta^k}{\delta t^k}(\frac{\frac{\delta \theta_j}{\delta t}}{\theta_j}))=\frac{\delta^{k+1}}{\delta t^{k+1}}(\frac{\frac{\delta \theta_j}{\delta u}}{\theta_j})$$ hold for $k=0,1,2,...$.}\\

 Indeed, we prove this lemma by recurrence after remarking that for $k=0$ we get $\frac{\delta}{\delta t}(\frac{\frac{\delta \theta_j}{\delta u}}{\theta_j})=\frac{\frac{\delta^2 \theta_j}{\delta u\delta t}}{\theta_j}-\frac{\delta \theta_j}{\delta u}\frac{\frac{\delta \theta_j}{\delta t}}{\theta^2_j}= \frac{\delta }{\delta u}(\frac{\frac{\delta \theta_j}{\delta t}}{\theta_j})$.\\
 
  By this lemma  $\frac{\delta }{\delta u}(\frac{\delta^k}{\delta t^k}(\frac{\frac{\delta \theta_j}{\delta t}}{\theta_j})) < (>) 0$ means that $\frac{\delta^k}{\delta t^k}(\frac{\frac{\delta \theta_j}{\delta t}}{\theta_j})$ is decreasing (increasing) as a function of $u$.\\ 

{\bf Proof of Theorem 5}\quad   
By Theorem 4, $\frac{\frac{\delta \theta_{j}(u, t)}{\delta u}}{\theta_j(u, t)}$ are completely monotonic with respect to $t > 0$ for $j=1,4$.Thus, $(-1)^k(\frac{\delta^k}{\delta t^k}\frac{\frac{\delta \theta_{j}(u, t)}{\delta u}}{\theta_j(u, t)}) =(-1)^k \frac{\delta }{\delta u}(\frac{\delta^{k-1}}{\delta t^{k-1}}(\frac{\frac{\delta \theta_j}{\delta t}}{\theta_j})) > 0$. That means $(-1)^{k} \frac{\delta^k}{\delta t^k}(\frac{\frac{\delta S_j}{\delta t}}{S_j}) >  0$ hold for $k=0,1,2,...$. While $-\frac{\theta'_j(u,i\pi t)}{\theta_j(u,i\pi t)}$ are completely monotonic with respect to $t > 0$ for $j=2,3$ means $(-1)^{k} \frac{\delta^k}{\delta t^k}(\frac{\frac{\delta S_j}{\delta t}}{S_j}) <  0$ for $k=0,1,2,...$. So, Theorem 5 is proved\\

\section{Further results}
In this part, we will deduce consequences of Theorems 4 and 5 and find again some known results.\\ 
At first, we prove the following which has been proved in [7], [10] \\

{\bf Corollary 5} \quad {\it For fixed $u,v$ such that $0\leq u < v < 1$, the functions $S_1(u,v;t)$ and $S_4(u,v;t)$ are positive and strictly increasing for $t > 0$ while $S_2(u,v;t)$ and $S_3(u,v;t)$ are positive decreasing for $t > 0$.}

\bigskip
 Indeed, since $$\frac{\delta}{\delta t}\frac{\frac{\delta \theta_j}{\delta u}}{\theta_j}=\frac{\delta}{\delta u}\frac{\frac{\delta \theta_j}{\delta t}}{\theta_j}$$ and $\frac{\frac{\delta \theta_j}{\delta u}}{\theta_j}$ is CM then $\frac{\delta}{\delta u}\frac{\frac{\delta \theta_j}{\delta t}}{\theta_j} < 0$ for $j=1,4$. This implies that  
 $\frac{\frac{\delta \theta_j}{\delta t}}{\theta_j}$ is decreasing as a function of $u$. Then for $j=1,4$ the quotient $S_j(u,v,t)$ is positive and strictly increasing for $t > 0$. While $-\frac{\frac{\delta \theta_j}{\delta u}}{\theta_j}$ is CM then $\frac{\delta}{\delta u}\frac{\frac{\delta \theta_j}{\delta t}}{\theta_j} > 0$ for $j=2,3$. This implies that  
 $\frac{\frac{\delta \theta_j}{\delta t}}{\theta_j}$ is increasing as a function of $u$. Then for $j=2,3$ the quotient  $S_j(u,v,t)$ is positive and strictly decreasing for $t > 0$. Corollary 5 is then proved.\\

\bigskip

 We may also deduce the following which extends Theorem 1.2 of [10]
  \bigskip
 
 {\bf Corollary 6}\quad {\it For any fixed $u, v$ such that $0<u<v<1$ the functions $S_j$ are such that \\
 - For any $t > 0$ the second derivatives $\frac{\delta^2}{\delta t^2}S_j(u,v,t) < 0 $ for $j=1,4$.\\
 - For any $t > 0$ the second derivatives $\frac{\delta^2}{\delta t^2}S_j(u,v,t) > 0 $ for $j=2,3$.}\\
 
 \bigskip

{\bf Remark} \quad We may prove Corollary 6 thank to {\it Maple} or {\it Mathematica}. It allows us to compute the derivative with respect to $u$ of $\frac{d}{dt}(\frac{\theta'_j}{\theta_j})$ for $j =1,2,3,4.$

$$\frac{d}{du}[\frac{d}{dt}(\frac{\theta'_3}{\theta_3})]= \sum_{k\geq 0} {\frac {\cos \left( \frac{\pi}{2} \,u \right)  \left( 1-
 \left( \coth \left( \frac{\pi}{2} \, \left( 2\,k+1 \right) t \right) 
 \right) ^{2} \right) {\pi }^{3} \left( 2\,k+1 \right) ^{2}\sin
 \left( \frac{\pi}{2} \,u \right) }{\cosh \left( \pi \, \left( 2\,k+1
 \right) t \right) +\cos \left( \pi \,u \right) }}-$$ $${\frac { \left( 
\cos \left( \frac{\pi}{2} \,u \right)  \right) ^{2} \left( 1- \left( \coth
 \left( \frac{\pi}{2} \, \left( 2\,k+1 \right) t \right)  \right) ^{2}
 \right) {\pi }^{3} \left( 2\,k+1 \right) ^{2}\sin \left( \pi \,u
 \right) }{ \left( \cosh \left( \pi \, \left( 2\,k+1 \right) t
 \right) +\cos \left( \pi \,u \right)  \right) ^{2}}}-$$ $$2\,{\frac {\cos
 \left( \frac{\pi}{2} \,u \right) \coth \left( \frac{\pi}{2} \, \left( 2\,k+1
 \right) t \right) {\pi }^{3} \left( 2\,k+1 \right) ^{2}\sinh \left( 
\pi \, \left( 2\,k+1 \right) t \right) \sin \left( \frac{\pi}{2} \,u
 \right) }{ \left( \cosh \left( \pi \, \left( 2\,k+1 \right) t
 \right) +\cos \left( \pi \,u \right)  \right) ^{2}}}+$$ $$4\,{\frac {
 \left( \cos \left( \frac{\pi}{2} \,u \right)  \right) ^{2}\coth \left( \frac{\pi}{2} \, \left( 2\,k+1 \right) t \right) {\pi }^{3} \left( 2\,k+1
 \right) ^{2}\sinh \left( \pi \, \left( 2\,k+1 \right) t \right) \sin
 \left( \pi \,u \right) }{ \left( \cosh \left( \pi \, \left( 2\,k+1
 \right) t \right) +\cos \left( \pi \,u \right)  \right) ^{3}}} 
$$

$$\frac{d}{du}[\frac{d}{dt}(\frac{\theta'_4}{\theta_4})]= \sum_{k\geq 0} -{\frac {\cos \left( \frac{\pi}{2} \,u \right)  \left( 1-
 \left( \coth \left( \frac{\pi}{2} \, \left( 2\,k+1 \right) t \right) 
 \right) ^{2} \right) {\pi }^{3} \left( 2\,k+1 \right) ^{2}\sin
 \left( \frac{\pi}{2} \,u \right) }{\cosh \left( \pi \, \left( 2\,k+1
 \right) t \right) -\cos \left( \pi \,u \right) }}+$$ $${\frac { \left( 
\sin \left( \frac{\pi}{2} \,u \right)  \right) ^{2} \left( 1- \left( \coth
 \left( \frac{\pi}{2} \, \left( 2\,k+1 \right) t \right)  \right) ^{2}
 \right) {\pi }^{3} \left( 2\,k+1 \right) ^{2}\sin \left( \pi \,u
 \right) }{ \left( \cosh \left( \pi \, \left( 2\,k+1 \right) t
 \right) -\cos \left( \pi \,u \right)  \right) ^{2}}}+$$ $$2\,{\frac {\cos
 \left( \frac{\pi}{2} \,u \right) \coth \left( \frac{\pi}{2} \, \left( 2\,k+1
 \right) t \right) {\pi }^{3} \left( 2\,k+1 \right) ^{2}\sinh \left( 
\pi \, \left( 2\,k+1 \right) t \right) \sin \left( \frac{\pi}{2} \,u
 \right) }{ \left( \cosh \left( \pi \, \left( 2\,k+1 \right) t
 \right) -\cos \left( \pi \,u \right)  \right) ^{2}}}-$$ $$4\,{\frac {
 \left( \sin \left( \frac{\pi}{2} \,u \right)  \right) ^{2}\coth \left( \frac{\pi}{2} \, \left( 2\,k+1 \right) t \right) {\pi }^{3} \left( 2\,k+1
 \right) ^{2}\sinh \left( \pi \, \left( 2\,k+1 \right) t \right) \sin
 \left( \pi \,u \right) }{ \left( \cosh \left( \pi \, \left( 2\,k+1
 \right) t \right) -\cos \left( \pi \,u \right)  \right) ^{3}}}  
$$ 

$$\frac{d}{du}[\frac{d}{dt}(\frac{\theta'_1}{\theta_1})]=\sum_{k\geq 0}-{\frac {\frac{i}{2} \left( \cos \left( \frac{\pi}{2} \, \left( u+\frac{i}{2}
t \right)  \right)  \right) ^{2}{\pi }^{3}\coth \left( \frac{\pi}{2} \,
 \left( 2\,k+1 \right) t \right)  \left( 2\,k+1 \right) }{\cosh
 \left( \pi \, \left( 2\,k+1 \right) t \right) -\cos \left( \pi \,
 \left( u+\frac{i}{2}t \right)  \right) }}+$$ $${\frac {i\sin \left( \frac{\pi}{2} \,
 \left( u+\frac{i}{2}t \right)  \right) \coth \left( \frac{\pi}{2} \, \left( 2\,
k+1 \right) t \right) {\pi }^{3} \left( 2\,k+1 \right) \cos \left( 1/2
\,\pi \, \left( u+\frac{i}{2}t \right)  \right) \sin \left( \pi \, \left( u
+\frac{i}{2}t \right)  \right) }{ \left( \cosh \left( \pi \, \left( 2\,k+1
 \right) t \right) -\cos \left( \pi \, \left( u+\frac{i}{2}t \right) 
 \right)  \right) ^{2}}}+$$ $${\frac {\frac{i}{2} \left( \sin \left( \frac{\pi}{2} \,
 \left( u+\frac{i}{2}t \right)  \right)  \right) ^{2}\coth \left( \frac{\pi}{2} 
\, \left( 2\,k+1 \right) t \right) {\pi }^{3} \left( 2\,k+1 \right) }{
\cosh \left( \pi \, \left( 2\,k+1 \right) t \right) -\cos \left( \pi 
\, \left( u+\frac{i}{2}t \right)  \right) }}-$$ $${\frac {\sin \left( \frac{\pi}{2} 
\, \left( u+\frac{i}{2}t \right)  \right)  \left( 1- \left( \coth \left( \frac{\pi}{2} \, \left( 2\,k+1 \right) t \right)  \right) ^{2} \right) {\pi }
^{3} \left( 2\,k+1 \right) ^{2}\cos \left( \frac{\pi}{2} \, \left( u+\frac{i}{2}
t \right)  \right) }{\cosh \left( \pi \, \left( 2\,k+1 \right) t
 \right) -\cos \left( \pi \, \left( u+\frac{i}{2}t \right)  \right) }}+$$ $${
\frac { \left( \sin \left( \frac{\pi}{2} \, \left( u+\frac{i}{2}t \right) 
 \right)  \right) ^{2} \left( 1- \left( \coth \left( \frac{\pi}{2} \,
 \left( 2\,k+1 \right) t \right)  \right) ^{2} \right) {\pi }^{3}
 \left( 2\,k+1 \right) ^{2}\sin \left( \pi \, \left( u+\frac{i}{2}t
 \right)  \right) }{ \left( \cosh \left( \pi \, \left( 2\,k+1 \right) 
t \right) -\cos \left( \pi \, \left( u+\frac{i}{2}t \right)  \right) 
 \right) ^{2}}}+$$ $$2\,{\frac {\sin \left( \frac{\pi}{2} \, \left( u+\frac{i}{2}t
 \right)  \right) \coth \left( \frac{\pi}{2} \, \left( 2\,k+1 \right) t
 \right) {\pi }^{2} \left( 2\,k+1 \right)  \left( \sinh \left( \pi \,
 \left( 2\,k+1 \right) t \right) \pi \, \left( 2\,k+1 \right) +\frac{i}{2}
\sin \left( \pi \, \left( u+\frac{i}{2}t \right)  \right) \pi  \right) \cos
 \left( \frac{\pi}{2} \, \left( u+\frac{i}{2}t \right)  \right) }{ \left( \cosh
 \left( \pi \, \left( 2\,k+1 \right) t \right) -\cos \left( \pi \,
 \left( u+\frac{i}{2}t \right)  \right)  \right) ^{2}}}-$$ $$4\,{\frac { \left( 
\sin \left( \frac{\pi}{2} \, \left( u+\frac{i}{2}t \right)  \right)  \right) ^{2
}\coth \left( \frac{\pi}{2} \, \left( 2\,k+1 \right) t \right) {\pi }^{2}
 \left( 2\,k+1 \right)  \left( \sinh \left( \pi \, \left( 2\,k+1
 \right) t \right) \pi \, \left( 2\,k+1 \right) +\frac{i}{2}\sin \left( 
\pi \, \left( u+\frac{i}{2}t \right)  \right) \pi  \right) \sin \left( \pi 
\, \left( u+\frac{i}{2}t \right)  \right) }{ \left( \cosh \left( \pi \,
 \left( 2\,k+1 \right) t \right) -\cos \left( \pi \, \left( u+\frac{i}{2}t
 \right)  \right)  \right) ^{3}}}+$$ $${\frac {i \left( \sin \left( \frac{\pi}{2} \, \left( u+\frac{i}{2}t \right)  \right)  \right) ^{2}\coth \left( \frac{\pi}{2} \, \left( 2\,k+1 \right) t \right) {\pi }^{3} \left( 2\,k+1
 \right) \cos \left( \pi \, \left( u+\frac{i}{2}t \right)  \right) }{
 \left( \cosh \left( \pi \, \left( 2\,k+1 \right) t \right) -\cos
 \left( \pi \, \left( u+\frac{i}{2}t \right)  \right)  \right) ^{2}}}
$$

$$\frac{d}{du}[\frac{d}{dt}(\frac{\theta'_2}{\theta_2})]=\sum_{k\geq 0}-{\frac {\frac{i}{2} \left( \sin \left( \frac{\pi}{2} \, \left( u+\frac{i}{2}
t \right)  \right)  \right) ^{2}\coth \left( \frac{\pi}{2} \, \left( 2\,k+1
 \right) t \right) {\pi }^{3} \left( 2\,k+1 \right) }{\cosh \left( 
\pi \, \left( 2\,k+1 \right) t \right) +\cos \left( \pi \, \left( u+1/
2\,it \right)  \right) }}+$$ $${\frac {i\sin \left( \frac{\pi}{2} \, \left( u+1/
2\,it \right)  \right) \coth \left( \frac{\pi}{2} \, \left( 2\,k+1 \right) 
t \right) {\pi }^{3} \left( 2\,k+1 \right) \cos \left( \frac{\pi}{2} \,
 \left( u+\frac{i}{2}t \right)  \right) \sin \left( \pi \, \left( u+\frac{i}{2}t
 \right)  \right) }{ \left( \cosh \left( \pi \, \left( 2\,k+1 \right) 
t \right) +\cos \left( \pi \, \left( u+\frac{i}{2}t \right)  \right) 
 \right) ^{2}}}+$$ $${\frac {\frac{i}{2} \left( \cos \left( \frac{\pi}{2} \, \left( u
+\frac{i}{2}t \right)  \right)  \right) ^{2}{\pi }^{3}\coth \left( 1/2\,
\pi \, \left( 2\,k+1 \right) t \right)  \left( 2\,k+1 \right) }{\cosh
 \left( \pi \, \left( 2\,k+1 \right) t \right) +\cos \left( \pi \,
 \left( u+\frac{i}{2}t \right)  \right) }}+$$ $${\frac {\sin \left( \frac{\pi}{2} \,
 \left( u+\frac{i}{2}t \right)  \right)  \left( 1- \left( \coth \left( \frac{\pi}{2} \, \left( 2\,k+1 \right) t \right)  \right) ^{2} \right) {\pi }^
{3} \left( 2\,k+1 \right) ^{2}\cos \left( \frac{\pi}{2} \, \left( u+\frac{i}{2}t
 \right)  \right) }{\cosh \left( \pi \, \left( 2\,k+1 \right) t
 \right) +\cos \left( \pi \, \left( u+\frac{i}{2}t \right)  \right) }}-$$ $${
\frac { \left( \cos \left( \frac{\pi}{2} \, \left( u+\frac{i}{2}t \right) 
 \right)  \right) ^{2} \left( 1- \left( \coth \left( \frac{\pi}{2} \,
 \left( 2\,k+1 \right) t \right)  \right) ^{2} \right) {\pi }^{3}
 \left( 2\,k+1 \right) ^{2}\sin \left( \pi \, \left( u+\frac{i}{2}t
 \right)  \right) }{ \left( \cosh \left( \pi \, \left( 2\,k+1 \right) 
t \right) +\cos \left( \pi \, \left( u+\frac{i}{2}t \right)  \right) 
 \right) ^{2}}}-$$ $$2\,{\frac {\cos \left( \frac{\pi}{2} \, \left( u+\frac{i}{2}t
 \right)  \right) \coth \left( \frac{\pi}{2} \, \left( 2\,k+1 \right) t
 \right) {\pi }^{2} \left( 2\,k+1 \right)  \left( \sinh \left( \pi \,
 \left( 2\,k+1 \right) t \right) \pi \, \left( 2\,k+1 \right) -\frac{i}{2}
\sin \left( \pi \, \left( u+\frac{i}{2}t \right)  \right) \pi  \right) \sin
 \left( \frac{\pi}{2} \, \left( u+\frac{i}{2}t \right)  \right) }{ \left( \cosh
 \left( \pi \, \left( 2\,k+1 \right) t \right) +\cos \left( \pi \,
 \left( u+\frac{i}{2}t \right)  \right)  \right) ^{2}}}+$$ $$4\,{\frac { \left( 
\cos \left( \frac{\pi}{2} \, \left( u+\frac{i}{2}t \right)  \right)  \right) ^{2
}\coth \left( \frac{\pi}{2} \, \left( 2\,k+1 \right) t \right) {\pi }^{2}
 \left( 2\,k+1 \right)  \left( \sinh \left( \pi \, \left( 2\,k+1
 \right) t \right) \pi \, \left( 2\,k+1 \right) -\frac{i}{2}\sin \left( 
\pi \, \left( u+\frac{i}{2}t \right)  \right) \pi  \right) \sin \left( \pi 
\, \left( u+\frac{i}{2}t \right)  \right) }{ \left( \cosh \left( \pi \,
 \left( 2\,k+1 \right) t \right) +\cos \left( \pi \, \left( u+\frac{i}{2}t
 \right)  \right)  \right) ^{3}}}-$$ $${\frac {i \left( \cos \left( \frac{\pi}{2} \, \left( u+\frac{i}{2}t \right)  \right)  \right) ^{2}\coth \left( \frac{\pi}{2} \, \left( 2\,k+1 \right) t \right) {\pi }^{3} \left( 2\,k+1
 \right) \cos \left( \pi \, \left( u+\frac{i}{2}t \right)  \right) }{
 \left( \cosh \left( \pi \, \left( 2\,k+1 \right) t \right) +\cos
 \left( \pi \, \left( u+\frac{i}{2}t \right)  \right)  \right) ^{2}}}
$$
Again thank to {\it Maple} or {\it Mathematica} we are able to verify that \  
 $$\frac{d}{du}[\frac{d}{dt}(\frac{\theta'_j}{\theta_j})], j=1,4 \ {\mbox is\ positive\ while}\quad \frac{d}{du}[\frac{d}{dt}(\frac{\theta'_2}{\theta_2})], j=2,3 \ is \ negative.$$  This means that \ 
 $$\frac{S'_j}{S_j}(u,v,t), j=1,4 \ {\mbox is\ positive\ while}\quad \frac{S'_j}{S_j}(u,v,t), j=2,3\ is\ negative.$$\\

\bigskip

{\bf Corollary 7}\quad {\it  For any fixed $u, v$ such that $0<u<v<1$ the function $S_j$ for $j=2,3$ is such that \ 
 $\frac{\delta^3}{\delta t^3} S_j(u,v;t)$ is non positive for $t \in ]0,\infty[$}\\
 
 Indeed,  
 $$\frac{\delta}{\delta t}\frac{\frac{\delta}{\delta t} S_j}{S_j}(u,v,t)=\frac{\delta}{\delta t}(\frac{\theta'_j}{\theta_j})(u, t)-\frac{d}{dt}(\frac{\theta'_j}{\theta_j})(v, t)$$ implies
 $$\frac{\delta^2}{\delta t^2}\frac{\frac{\delta}{\delta t} S_j}{S_j}(u,v,t)=\frac{\delta^2}{\delta t^2}(\frac{\theta'_j}{\theta_j})(u, t)-\frac{\delta^2}{\delta t^2}(\frac{\theta'_j}{\theta_j})(v, t).$$
 Since by Lemma 2 $$\frac{\delta^3}{\delta t^3}\frac{\frac{\delta \theta'_j}{\delta u}}{\theta_j}=\frac{\delta}{\delta u}\frac{\delta^2}{\delta t^2}(\frac{\theta'_j}{\theta_j})$$ and $-\frac{\frac{\delta \theta_j}{\delta u}}{\theta_j}$ is CM then $\frac{\delta}{\delta u}\frac{\delta^2}{\delta t^2}(\frac{\theta'_j}{\theta_j}) < 0$ for $j=2,3$. This implies that $\frac{\delta^2}{\delta t^2}(\frac{\theta'_j}{\theta_j})$ is decreasing as function of $u$. Thus,   
 $$\frac{\delta^2}{\delta t^2}\frac{\frac{\delta}{\delta t} S_j}{S_j}=\frac{\frac{\delta^3}{\delta t^3}S}{S}-3\frac{\frac{\delta}{\delta t}S'\frac{\delta^2}{\delta t^2}S}{S^2}+2\frac{(\frac{\delta}{\delta t}S)^2}{S^3} < 0.$$ By Corollaries 5 and 6, we have  $\frac{\delta}{\delta t}S < 0$ and $\frac{\delta^2}{\delta t^2}S > 0$ for $j=2,3.$ Then we get $\frac{\delta^3}{\delta t^3}S < 0$.\\

The following only concerns the quotient $S_4(u,v,t)$\\
 
 {\bf Corollary 8}\quad {\it  For any fixed $u, v$ such that $0<u<v<1$ the function $S_4(u,v,t)$ is such that \ 
 $\frac{d^2}{dt^2}(\frac{\frac{\delta}{\delta t}S_4}{S_4})$ is non negative for $t \in ]0,\infty[$}
 
 \bigskip
 Notice that since $$\frac{\frac{\delta}{\delta t} S_4}{S_4}(u,v,t) = \frac{\theta'_4}{\theta_4}(u, t) -\frac{\theta'_4}{\theta_4}(v,i \pi t)$$ then $$\frac{\delta^2}{\delta t^2}\frac{\frac{\delta}{\delta t} S_j}{S_4}(u,v,t)=\frac{\delta^2}{\delta t^2}(\frac{\theta'_4}{\theta_4})(u, t)-\frac{\delta^2}{\delta t^2}(\frac{\theta'_4}{\theta_4})(v, t).$$
 This Corollary is a direct consequence of Theorem 5. But we will use {\it Maple} or {\it Mathematica}.
 Compute for that next derivatives
$$[\frac{d^2}{dt^2}(\frac{\theta'_4}{\theta_4})]= \sum_{k\geq 0} {\frac { \left( \sin \left( 1/2
\,\pi \,u \right)  \right) ^{2}\coth \left( \frac{\pi}{2} \, \left( 2\,k+1
 \right) t \right)  \left( 1- \left( \coth \left( \frac{\pi}{2} \, \left( 2
\,k+1 \right) t \right)  \right) ^{2} \right) {\pi }^{3} \left( 2\,k+1
 \right) ^{3}}{\cosh \left( \pi \, \left( 2\,k+1 \right) t \right) -
\cos \left( \pi \,u \right) }}+$$ $$2\,{\frac { \left( \sin \left( 1/2\,
\pi \,u \right)  \right) ^{2} \left( 1- \left( \coth \left( \frac{\pi}{2} 
\, \left( 2\,k+1 \right) t \right)  \right) ^{2} \right) {\pi }^{3}
 \left( 2\,k+1 \right) ^{3}\sinh \left( \pi \, \left( 2\,k+1 \right) t
 \right) }{ \left( \cosh \left( \pi \, \left( 2\,k+1 \right) t
 \right) -\cos \left( \pi \,u \right)  \right) ^{2}}}-$$ $$ 4\,{\frac {
 \left( \sin \left( \frac{\pi}{2} \,u \right)  \right) ^{2}\coth \left( 1/2
\,\pi \, \left( 2\,k+1 \right) t \right) {\pi }^{3} \left( 2\,k+1
 \right) ^{3} \left( \sinh \left( \pi \, \left( 2\,k+1 \right) t
 \right)  \right) ^{2}}{ \left( \cosh \left( \pi \, \left( 2\,k+1
 \right) t \right) -\cos \left( \pi \,u \right)  \right) ^{3}}}+$$ $$2\,{
\frac { \left( \sin \left( \frac{\pi}{2} \,u \right)  \right) ^{2}\coth
 \left( \frac{\pi}{2} \, \left( 2\,k+1 \right) t \right) {\pi }^{3} \left( 
2\,k+1 \right) ^{3}\cosh \left( \pi \, \left( 2\,k+1 \right) t
 \right) }{ \left( \cosh \left( \pi \, \left( 2\,k+1 \right) t
 \right) -\cos \left( \pi \,u \right)  \right) ^{2}}} 
$$
as well as the derivative with respect to $u$ gives $$
\frac{d}{du}[\frac{d^2}{dt^2}(\frac{\theta'_4}{\theta_4})]= \sum_{k\geq 0} {\frac {\sin \left( \frac{\pi}{2} \,u \right) \coth \left( 
\frac{\pi}{2} \, \left( 2\,k+1 \right) t \right)  \left( 1- \left( \coth
 \left( \frac{\pi}{2} \, \left( 2\,k+1 \right) t \right)  \right) ^{2}
 \right) {\pi }^{4} \left( 2\,k+1 \right) ^{3}\cos \left( \frac{\pi}{2} \,u
 \right) }{\cosh \left( \pi \, \left( 2\,k+1 \right) t \right) -\cos
 \left( \pi \,u \right) }}-$$ $${\frac { \left( \sin \left( \frac{\pi}{2} \,u
 \right)  \right) ^{2}\coth \left( \frac{\pi}{2} \, \left( 2\,k+1 \right) t
 \right)  \left( 1- \left( \coth \left( \frac{\pi}{2} \, \left( 2\,k+1
 \right) t \right)  \right) ^{2} \right) {\pi }^{4} \left( 2\,k+1
 \right) ^{3}\sin \left( \pi \,u \right) }{ \left( \cosh \left( \pi \,
 \left( 2\,k+1 \right) t \right) -\cos \left( \pi \,u \right) 
 \right) ^{2}}}+$$ $$2\,{\frac {\sin \left( \frac{\pi}{2} \,u \right)  \left( 1-
 \left( \coth \left( \frac{\pi}{2} \, \left( 2\,k+1 \right) t \right) 
 \right) ^{2} \right) {\pi }^{4} \left( 2\,k+1 \right) ^{3}\sinh
 \left( \pi \, \left( 2\,k+1 \right) t \right) \cos \left( \frac{\pi}{2} \,
u \right) }{ \left( \cosh \left( \pi \, \left( 2\,k+1 \right) t
 \right) -\cos \left( \pi \,u \right)  \right) ^{2}}}-$$ $$4\,{\frac {
 \left( \sin \left( \frac{\pi}{2} \,u \right)  \right) ^{2} \left( 1-
 \left( \coth \left( \frac{\pi}{2} \, \left( 2\,k+1 \right) t \right) 
 \right) ^{2} \right) {\pi }^{4} \left( 2\,k+1 \right) ^{3}\sinh
 \left( \pi \, \left( 2\,k+1 \right) t \right) \sin \left( \pi \,u
 \right) }{ \left( \cosh \left( \pi \, \left( 2\,k+1 \right) t
 \right) -\cos \left( \pi \,u \right)  \right) ^{3}}}-$$ $$4\,{\frac {\sin
 \left( \frac{\pi}{2} \,u \right) \coth \left( \frac{\pi}{2} \, \left( 2\,k+1
 \right) t \right) {\pi }^{4} \left( 2\,k+1 \right) ^{3} \left( \sinh
 \left( \pi \, \left( 2\,k+1 \right) t \right)  \right) ^{2}\cos
 \left( \frac{\pi}{2} \,u \right) }{ \left( \cosh \left( \pi \, \left( 2\,k
+1 \right) t \right) -\cos \left( \pi \,u \right)  \right) ^{3}}}+$$ $$12\,
{\frac { \left( \sin \left( \frac{\pi}{2} \,u \right)  \right) ^{2}\coth
 \left( \frac{\pi}{2} \, \left( 2\,k+1 \right) t \right) {\pi }^{4} \left( 
2\,k+1 \right) ^{3} \left( \sinh \left( \pi \, \left( 2\,k+1 \right) t
 \right)  \right) ^{2}\sin \left( \pi \,u \right) }{ \left( \cosh
 \left( \pi \, \left( 2\,k+1 \right) t \right) -\cos \left( \pi \,u
 \right)  \right) ^{4}}}+$$ $$2\,{\frac {\sin \left( \frac{\pi}{2} \,u \right) 
\coth \left( \frac{\pi}{2} \, \left( 2\,k+1 \right) t \right) {\pi }^{4}
 \left( 2\,k+1 \right) ^{3}\cosh \left( \pi \, \left( 2\,k+1 \right) t
 \right) \cos \left( \frac{\pi}{2} \,u \right) }{ \left( \cosh \left( \pi 
\, \left( 2\,k+1 \right) t \right) -\cos \left( \pi \,u \right) 
 \right) ^{2}}}-$$ $$4\,{\frac { \left( \sin \left( \frac{\pi}{2} \,u \right) 
 \right) ^{2}\coth \left( \frac{\pi}{2} \, \left( 2\,k+1 \right) t \right) 
{\pi }^{4} \left( 2\,k+1 \right) ^{3}\cosh \left( \pi \, \left( 2\,k+1
 \right) t \right) \sin \left( \pi \,u \right) }{ \left( \cosh \left( 
\pi \, \left( 2\,k+1 \right) t \right) -\cos \left( \pi \,u \right) 
 \right) ^{3}}} 
$$
Thanks to {\it Maple} or {\it Mathematica} we verify that for any fixed $u$ and $t > 0$ that $\frac{d}{du}[\frac{d^2}{dt^2}(\frac{\theta'_4}{\theta_4})]$ is negative. This means that $\frac{d^2}{dt^2}(\frac{\theta'_4}{\theta_4})$ decreases as a function of $u$. Then $$\frac{d^2}{dt^2}(\frac{S'_4}{S_4})=\frac{d^2}{dt^2}(\frac{\theta'_4}{\theta_4})-\frac{d^2}{dt^2}(\frac{\theta'_4}{\theta_4})$$ is non negative. 

\newpage

{\bf References}\\

[1] P. Appell, E. Lacour, \quad {\it Fonctions elliptiques et applications}\quad Gauthiers-Villard ed., Paris (1922).\\

[2] H. Alzer and C. Berg, \quad {\it Some classes of completely monotonic functions,} Ann. Acad. Scient. Fennicae 27 (2002), 445-460.\\

[3] \ W. Cheney, W. Light, \quad {\it A course in Approximation theory}, Grad. texts in Math., vol 101, A.M.S., Providence, 2009.\\

[4] \ A.R. Chouikha, \quad {\it On Properties of Elliptic Jacobi Functions and Applications}, J. of Nonlin. Math. Physics Vol. 12-2, (2005), 162-169.\\

[5]\ A.R. Chouikha \quad \textit{Expansions of Theta Functions and Applications} \ ArXiv, math/0112137,  http://front.math.ucdavis.edu/0112.5137, (2011).\\

[6]\ A. Dixit, A. Solynin \quad {\it Monotonicity of quotients of theta functions related to an extremal problem on harmonic measure}, J. Math. Anal. Appl., 336 (2007), 1042-1053.\\

[7] A. Dixit, A. Roy, A. Zaharescu \quad {\it Convexity of quotients of theta functions}, \  
J. Math. Anal. Appl., 386 (2012), 319-331.\\

[8] A.Erdélyi, W.Magnus, F.Oberhettinger,  F.Tricomi, \quad {\it  Higher transcendental
functions}\quad  Vol. II. Based on notes left by H. Bateman. Robert E. Krieger Publish. Co.,
Inc., Melbourne, Fla., (1981).\\

[9] Chun-Fu Wei1, Bai-Ni Guo, \quad  {\it Complete Monotonicity of Functions Connected with
the Exponential Function and Derivatives} \quad Abstr. and Appl. Analysis, vol 2014, article ID 851213.\\

[10] A.Yu. Solynin \quad {\it Harmonic measure of radial line segments and symmetrization} \quad Math. Sb. 189 (11-12), (1998), 1701-1718.\\

[11] E.T. Whittaker,G.N. Watson \quad {\it A course of Modern Analysis}\\ Cambridge (1963).\\ 

[12] D.V. Widder, \ {\it The Laplace Transform}, Princeton Univ. Press, Princeton, NJ, 1941.

\end{document}